\documentclass[12pt]{article}
\usepackage{graphicx}
\usepackage[cmtip,arrow]{xy}
\usepackage{amsmath,amssymb,pb-diagram,pb-xy}
\def\today{\number\day .\number\month .\number\year}
\usepackage[applemac]{inputenc}

\def \bs{\backslash}
\def \C{{\mathbb C}}

\def \CF{{\cal F}}

\def \df{\ \begin{array}{c} _{\rm def}\\ ^{\displaystyle =}\end{array}\ }

\def \Ga{\Gamma}
\def \ga{\gamma}
\def \H{{\mathbb H}}

\def \Im{{\rm Im}}

\def \la{\lambda}

\def \prf{{\bf Proof: }}
\def \PSL{{\rm PSL}}
\def \Q{{\mathbb Q}}
\def \qed{\ifmmode\eqno \square 
		\else\noproof\vskip 12pt plus 3pt minus 9pt \fi}
\def \noproof{{\unskip\nobreak\hfill\penalty50\hskip2em\hbox{}%
     \nobreak\hfill $\square$\parfillskip=0pt%
     \finalhyphendemerits=0\par}}
\def \R{{\mathbb R}}

\def \setminus{\smallsetminus}
\def \SL{{\rm {SL}}}

\def \({\left(}
\def \){\right)}
\def \={{\ =\ }}

\newcommand{\norm}[1]{\left|\hspace{-1pt}\left| #1\right|\hspace{-1pt}\right|}

\renewcommand{\sp}[1]{\left\langle #1\right\rangle}
\newcommand{\ol}[1]{\overline{#1}}

\newtheorem{theorem}{Theorem}[section]

\newtheorem{lemma}[theorem]{Lemma}

\begin{document}

\pagestyle{myheadings} \markright{GROWTH of MODULAR SYMBOLS}

\title{On the growth of modular symbols}
\author{Anton Deitmar}
\date{}
\maketitle


$$ $$


\section*{Introduction}

Let $f$ be a holomorphic cusp form of weight two for the group $\Ga=\Ga_0(N)$.
Then $f(z)dz$ is a $\Ga$-invariant holomorphic differential on the upper half plane $\H$ in $\C$.
For $\ga\in\Ga$ define the modular symbol
$$
\sp{\ga,f}\df -2\pi i\int_{z_0}^{\ga z_0}f(z)\,dz,
$$
which is independent of the choice of the point $z_0\in\H\cup\Q\cup\{ i\infty\}$. 
The modular symbol for fixed $f$ is a group homomorphism from $\Ga$ to the additive group of complex numbers.
For $\ga=\left(\begin{array}{cc}a & b \\c & d\end{array}\right)\in\Ga$ let $\norm\ga=\max(|a|,|b|,|c|,|d|)$.
In this paper we show that for given $f$ the modular symbol as  logarithmic growth, i.e.,
$$
|\sp{\ga,f}|\ \le\ A\log\norm\ga +B
$$
for some $A,B\ge 0$.
In \cite{Gold1,Gold2}, D. Goldfeld conjectured that the modular symbol has moderate growth if also $f$ is allowed to vary among the normalized newforms.
The result of the present paper reduces Goldfeld's conjecture to a statement on the growth of the modular symbol on a set of generators of the group $\Ga_0(N)$.

Note that, as $f$ is a cusp form, the modular symbol vanishes on parabolic elements, that is, $\sp{p,f}=0$ for every parabolic element $p$ of $\Ga$.

\section{Growth of an additive homomorphism}
Let $\SL_2(\R)$ denote the group of real $2\times 2$ matrices of determinant one.
Let $G$ be the group $\PSL_2(\R)=\SL_2(\R)/\pm 1$.

For $x=\left(\begin{array}{cc}a & b \\c & d\end{array}\right)\in G$ let $\norm x=\max(|a|,|b|,|c|,|d|)$.
Let $M$ be a subset of $G$.
A function $f:M\to \C$ is said to be of \emph{logarithmic growth}, if there are constants $A,B\ge 0$ such that
$$
|f(x)|\ \le\ A\log\norm x +B
$$
holds for every $x\in M$.

\begin{theorem}
Let $\Ga$ be a lattice in $\PSL_2(\R)$ and let $\psi :\Ga\to\C$ be a group homomorphism with $\psi(p)=0$ for every parabolic element of $\Ga$.
Then $\psi$ is of logarithmic  growth.
\end{theorem}

\prf
Let $\H$ denote the upper half plane in $\C$ equipped with the hyperbolic metric $ds^2=\frac{dx^2+dy^2}{y^2}$.
Choose $z_0\in\H$ which is not a fixed point of an elliptic element of $\Ga$, and let $\CF$ denote the corresponding \emph{Dirichlet fundamental domain}, also called the \emph{Dirichlet polygon} \cite{Beard}, as it is a hyperbolic polygon with finitely many sides.
It is defined as
$$
\CF\=\left\{ z\in\H : d(z,z_0)<d(\ga z,z_0)\ \forall\ga\in\Ga\setminus\{ 1\} \right\},
$$
where $d(z,w)$ denotes the hyperbolic distance of two points in the upper half plane $\H$.
Let $S=S^{-1}$ be a finite set of generators of the group $\Ga$.
For $\ga\in\Ga$ we write $\ga=s_1\cdots s_n$ as a shortest word, so $n=l_S(\ga)$, the \emph{word length} of $\ga$ with respect to $S$.
Then,
$$
|\psi(\ga)|\=|\psi(s_1\cdots s_n)|\=\left|\sum_{j=1}^n\psi(s_j)\right| \le\ C_S\,l_S(\ga),
$$
where $C_S\ge 0$ is the maximum of the values $|\psi(s)|$ for $s\in S$.

We first consider the case of a uniform lattice $\Ga$, i.e., the quotient $\Ga\bs G$ is compact.
Then the closure $\bar\CF$ of $\CF$ in $\H$ is compact.
According to Theorem IV 23 of \cite{Harpe}, for every $z\in\H$ there exists $\la \ge 1$, $C\ge 0$ with $l_S(\ga)\le \la d(\ga z,z)+C$ for every $\ga\in\Ga$.
For $z=i$ this implies
$$
|\psi(\ga)|\ \le\ C_S\la d(\ga i,i)+C_SC.
$$
By Theorem 7.2.1 of \cite{Beard} one has for $z,w\in\H$,
$$
d(z,w)\=\log\( \frac{|z-\bar w|+|z-w|}{|z-\bar w|-|z-w|}\).
$$
Let $\ga=\left(\begin{array}{cc}a & b \\c & d\end{array}\right)\in\Ga$ we get
\begin{eqnarray*}
d(\ga i,i) &=& \log\( \frac{\left|\frac{ai+b}{ci+d}+i\right|+
\left|\frac{ai+b}{ci+d}-i\right|}
{\left|\frac{ai+b}{ci+d}+i\right|
-\left|\frac{ai+b}{ci+d}-i\right|}\) \\
&=& \log\( \frac{|ai+b-c+di|+|ai+b+c-di|}{|ai+b-c+di|-|ai+b+c-di|}\)\\
&=& \log\( \frac{\sqrt{(b-c)^2+(a+d)^2}+\sqrt{(b+c)^2+(a-d)^2}}{\sqrt{(b-c)^2+(a+d)^2}-\sqrt{(b+c)^2+(a-d)^2}}\).
\end{eqnarray*}
Now $(b-c)^2+(a+d)^2-((b+c)^1+(a-d)^2)=-4bc+4ad=4$, as the determinant of $\ga$ is 1.
Therefore,
\begin{eqnarray*}
d(\ga i,i) &=& \log\( \( \sqrt{(b-c)^2+(a+d)^2}+\sqrt{(b+c)^2+(a-d)^2}\)^2\)-2\log 2\\
&\le& 2\log\(\sqrt{8\norm\ga^2}+\sqrt{8\norm\ga^2}\) -2\log 2\\
&=& 2\log\norm\ga +3\log 2.
\end{eqnarray*}
This gives the Theorem in the case of $\Ga$ being a uniform lattice.

Next assume that $\Ga$ is not uniform.
Then the fundamental domain $\CF$ has cusps.
We assume that $z_0$ is chosen in a way that no two cusps of $\CF$ are equivalent under $\Ga$.
So the cusps of $\CF$ are a set of representatives of the set of cusps of $\Ga$ in the boundary $\partial\H$ of $\H$ module $\Ga$-equivalence.
For each cusp $c\in\partial\H$ fix some $\sigma_c\in G$ with $c=\sigma_c\infty$.
We do so in a $\Ga$-compatible way, i.e., for $\ga\in\Ga$ we suppose that $\sigma_{\ga c}\sigma_c^{-1}$ lies in $\Ga$.

Let $T>1$ and set
$$
\CF_T\=\{ z\in\CF : \Im(\sigma_c^{-1}z)\le T\ \mbox{for every cusp }c\}.
$$
We choose $T$ so large that $\CF_T$ equals $\CF$ minus cusp sections.
Let $\H_T$ be the union of all sets $\ol{\ga\CF_T}$ where $\ga$ ranges over $\Ga$.
Then $\H_T$ equals $\H$ minus a countable number of open horoballs.
Thus $\H_T$ is a Riemannian manifold with boundary.
Let $d_T$ denote the distance function on $\H_T$.
Note that if $z,w\in\H_T$, and the geodesic in $\H$ joining them lies completely in $\H_T$, then $d_T(z,w)=d(z,w)$.

For a cusp $c$ let 
$$
H_{c,T}\=\sigma_c(\{\Im(z)>T\})
$$
be the $T$-horoball attached to $c$.
Increasing $T$ if necessary, we can make sure, that the geodesic $\ol{z_0,\ga z_0}$ in $\H$ is disjoint to $\CF\cap H_{c,T}$ for every cusp $c$ of $\CF$ and every $\ga\in\Ga$.
Note that for every cusp $c$ of $\Ga$, which is not a cusp of $\CF$, the intersection $\CF\cap H_{c,T}$ is empty.

Let $\ga\in\Ga$ and suppose that the geodesic $\ol{z_0,\ga z_0}$ in $\H$ does not completely lie in $\H_T$.
Then this geodesic meets some horoball $H_{c,T}$.
After applying $\sigma_c^{-1}$, on can assume $c=\infty$.
Then there is a generator $p_c$ of the stabilizer group $\Ga_c$ of the cusp $c$, such that
the distance $d(z_0,p_c\ga z_0)$ is then strictly less than $d(z_0,\ga z_0)$.
From this it follows that if $d(z_0,\ga z_0)\le d(z_0,p\ga z_0)$ for every parabolic $p\in\Ga$, then the geodesic 
$\ol{z_0,\ga z_0}$ lies completely in the set $\H_T$.

By Theorem IV 23 in \cite{Harpe}, there are $\la\ge 1$ and $C\ge 0$ with $l_S(\ga)\le \la d_T(z_0,\ga z_0)+C$ for every $\ga\in\Ga$.

\begin{lemma}
Let $\ga\in\Ga\setminus\{ 1\}$ be given.
There are parabolic elements $p_1,\dots p_n$ of $\Ga$ such that with
$\ga_s=p_n\cdots p_1\ga$ one has $d(z_0,\ga_s z_0)\le d(z_0,\ga z_0)$ and the geodesic $\ol{z_0,\ga_s z_0}$ lies in $\H_T$.
\end{lemma}

\prf
Assume first that for every parabolic element $p\in\Ga$ the distance $d(z_0,\ga z_0)$ is less than or equal to $d(z_0,p\ga z_0)$.
Then the geodesic $\ol{z_0,\ga z_0}$ lies in $\H_T$.
We set $\ga_s=\ga$ and we are done.

Now if there exists a parabolic $p_1\in\Ga$ such that 
$d(z_0,p_1\ga z_0)<d(z_0,\ga z_0)$ then replace $\ga$ with $p_1\ga$.
After that, either the condition above is satisfied or we find a parabolic $p_2$ such that $d(z_0,p_2p_1\ga z_0)<d(z_0,p_1\ga z_)$.
Iteration yields a sequence $p_1,p_2,\dots\in\Ga$.
This process terminates, as for a given radius $r$ there are only finitely many $\Ga$-conjugates of $z_0$ in distance $\le r$.
\qed

To finish the proof of the theorem let $\ga\in\Ga$ and consider $\ga_s$ as in the lemma.
Then $|\psi(\ga)|\=|\psi(\ga_s)|\le C_S l_S(\ga_s)\le C_S\la d_T(z_0,\ga_s z_0)+C_SC$ and 
\begin{eqnarray*}
d_T(z_0,\ga_s z_0)&=& d(z_0,\ga_s z_0)\\
&\le & d(z_0,\ga z_0)\\
&\le& d(z_0,i)+d(i,\ga i)+d(\ga i,\ga z_0)\\
&\le& 2\log\norm\ga +3\log 2+ 2d(z_0,i).
\end{eqnarray*}
\qed

The value of these results with respect to the Goldfeld conjecture hinges on the control over the constants $C,C_S,\la$ as the group $\Ga$ shrinks.
The constant $C_S$ depends on the group $\Ga$ and on the homomorphism $\psi$, i.e., if $\psi$ is a modular symbol, on the cusp form $f$.
The constants $C$ and $\la$, however, do not depend on $\psi$, therefore they are easier to control.
The following explicit estimate might be useful.

Let $R$ be the diameter of $\CF_T$ and let $B$ be the closed ball in $\H_T$ around $z_0$ of radius $R$.
Let $S=\{ s\in\Ga : sB\cap B\ne \emptyset\}$.
Then $S=S^{-1}$ is a finite set of generators of $\Ga$.
Let 
$$
r\= \inf\{ d(B,\ga B): \ga\in\Ga\setminus S\}.
$$

\begin{lemma}
The number $r$ is $>0$ and
$$
|\psi(\ga)|\ \le\ C_S\( \frac{2}{r}\log\norm\ga + \frac{3\log 2}{r}+1\).
$$
\end{lemma}

\prf
Let $r_T=\inf\{ d_T(B,\ga B): \ga\in\Ga\setminus S\}$.
The proof of Theorem IV 23 of \cite{Harpe} together with the proof of our theorem yields
$$
|\psi(\ga)|\ \le\ \frac{C_S}{r_T}(2\log\norm\ga +3\log 2+2d(z_0,i))+C_S
$$
Now $r\ge r_T$ and the argument for $r_T>0$ also implies $r>0$.
Further, varying $z_0$ the distance $d(z_0,i)$ can be chosen arbitrarily small.
\qed

{\small Mathematisches Institut\\
Auf der Morgenstelle 10\\
72076 T\"ubingen\\
Germany\\
\tt deitmar@uni-tuebingen.de}

\today


\newpage
\begin{thebibliography}{XXX}

\bibitem{Beard}
\bf Beardon, A.F.:
\it The geometry of discrete groups.
\rm Graduate Texts in
Mathematics, 91. Springer-Verlag, New York, 1983.

\bibitem{Gold1}
\bf Goldfeld, D.:
\it Modular elliptic curves and Diophantine problems.
\rm Number theory (Banff, AB, 1988), 157--175, de Gruyter, Berlin, 1990.

\bibitem{Gold2}
\bf Goldfeld, D.:
\it Modular forms, elliptic curves and the $ABC$-conjecture.
\rm A panorama of number theory or the view from Baker's garden (Zürich, 1999), 128-147, Cambridge Univ. Press, Cambridge, 2002.

\bibitem{Harpe}
\bf de la Harpe, P.:
\it Topics in geometric group theory.
\rm Chicago Lectures in Mathematics. University of Chicago Press, Chicago, IL, 2000.

\bibitem{Shim}
\bf Shimura, G.:
\it On the factors of the jacobian variety of a modular function field.
\rm J. Math. Soc. Japan 25, 523-544 (1973).

\end{thebibliography}
\end{document}